# Entropy Rate of Thermal Diffusion

**John Laurence Haller Jr.**
Predictive Analytics, CCC Information Services, Chicago, USA
jlhaller@gmail.com

**Abstract**

The thermal diffusion of a free particle is a random process and generates entropy at a rate equal to twice the particle's temperature, $R = 2k_B T/\hbar$ (in natural units of information per second). The rate is calculated using a Gaussian process with a variance of $(\Delta x_0 + \Delta p \cdot t/m)^2$ which is a combination of quantum and classical diffusion. The solution to the quantum diffusion of a free particle is derived from the equation for kinetic energy and its associated imaginary diffusion constant; a real diffusion constant (representing classical diffusion) is shown to be $D = \hbar/(2m)$. We find the entropy of the initial state is one natural unit, which is the same amount of entropy the process generates after the de-coherence time, $\tau = \hbar/2k_B T$.

**1. Primary Finding:**

When a free particle is at a non-zero temperature, it is composed of a spectrum of frequencies that evolve at different rates which causes the probability distribution of where one can find the particle to spread. We will show that the entropy rate, associated with the probability distribution diffusing, is equal to twice the particle's temperature.

$$R = 2k_B T/\hbar \quad (1)$$

The rate, R, is calculated below using the natural logarithm, and thus the units for the rate are natural units of information per second, when the temperature (T) is expressed in degrees Kelvin, Boltzmann's constant ($k_B$) is expressed in Joules per Kelvin, and hbar ($\hbar$) is Planck's constant divided by $2\pi$ in Joule-seconds.

This equation tells us the minimum amount of information we need, each second, in order to track a diffusing free particle to the highest precision that nature requires. By quantifying the amount of information needed to follow a free particle for a certain time, and showing it is finite, we are able to guarantee that a computer (or other discrete state-space machine with finite memory) can store a particle's initial state and trajectory.

What is unique about this result is that there is no dependence on the mass of the particle or any other variable except the temperature.

**2. Assumptions:**

We prove this primary result by making the following three assumptions:
1) The continuous diffusion of a free particle can be modeled as a discrete process with a time step $dt$ that is much smaller than the de-coherence time $\tau$, $dt \ll \tau = \hbar/(2k_B T)$, where T is the temperature.
2) Knowing the particle's location at time step n+1 allows one to determine the location of the particle at the previous time step n; i.e., conditional entropy is zero, $h(X_n|X_{n+1}) = 0$ where $X_n$ is the random variable that represents where the particle can be found at time step n.
3) At each time step the minimum uncertainty wave-packet is localized around its new location, and thus the conditional entropy of the $n+1$ step, given all previous steps, is the same as the conditional entropy of the 1$^{st}$ step given its initial state, $h(X_{n+1}|X_n, X_{n-1}, \ldots, X_1, X_0) = h(X_1|X_0)$.

These three assumptions taken together are reasonable and give insight into the behavior of the system.

Assumption 1 is aided by the analysis found in [1] which shows the time step of discrete diffusion is $dt = \hbar/(2mc^2)$ and thus for non-relativistic particles, the assumption holds.

Assumption 2 says that there is no entropy beyond the minimum uncertainty wave-packet after a measurement of the particle's location was made.

Assumption 3 says that the vacuum localizes the diffusing particle up to the minimum uncertainty wave packet at each step in the process. Even though an un-disrupted particle's wave packet will spread, at the vacuum level the particle is re-initialized at each step, like in quantum nondemolition measurements [2].

**3. Setup:**

At $t = 0$, a free particle in vacuum is initialized into a minimum uncertainty Gaussian wave-packet with a spatial variance equal to $(\Delta x_0)^2$. As time increases so does its variance and thus its entropy.

To calculate the entropy rate of this process, it is helpful to think of time as occurring in discrete units of a small size dt (assumption 1).

We can look at a Venn diagram of this process, figure (1). $X_{g,0}$ (or X0 in the figure) is a random variable, drawn from $g(x, 0)$, that describes the location of where the particle can be found at time $t = 0$. $X_{g,1}$ (X1) is a random variable, drawn from $g(x, dt)$, that describes the location of where the particle can be found at time $t = dt$. $X_{g,2}$ (X2) is drawn from $g(x, 2 \cdot dt)$ and so on up to $X_{g,n}$ which is drawn from $g(x, t)$ where $t = n \cdot dt$.

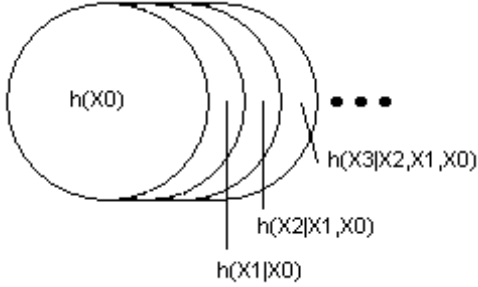

**Figure (1) – Venn diagram of the conditional entropies of the diffusion process**

As hinted to in the diagram (but explicitly stated here as assumption 2 and assumption 3), we will assume that the conditional entropy of each step is constant; $h(X_{n+1}|X_n, X_{n-1}, \dots, X_1, X_0) = h(X_1|X_0)$ and $h(X_n|X_{n+1}) = 0$, where h is the differential entropy $h = -\int g \cdot \log(g) \, dx$ and where $g(x,t)$ is the distribution which determines $X_n$. This essentially means that knowing the location of the particle at any time allows one to calculate where it was in the previous time step and that the minimum uncertainty wave-packet maintains its coherence as its first moment (or average value) diffuses via a process with a variance as given by equation 2.

In section 5, we show that as time increases, a free particle diffuses such that the variance of where the particle can be found (if localized) is $(\Delta x)^2$.

$$(\Delta x)^2 = \left(\Delta x_0 + \frac{\Delta p \cdot t}{m}\right)^2 \quad (2)$$

Thus $X_{g,n}$ (or simply $X_n$) is a Gaussian random variable with variance $(\Delta x)^2 = (\Delta x_0 + \Delta p \cdot n \cdot dt/m)^2$.

## 4. Entropy Rate:

We can calculate the entropy rate of this process using the definition of the entropy rate. We will use the entropy rate, R, as calculated by taking the limit as the number of steps goes to infinity of the conditional entropy of the last step given all previous steps divided by the time step [3].

$$R \equiv \lim_{n \to \infty} \frac{h(X_{n+1}|X_n, X_{n-1}, \dots, X_1, X_0)}{dt} \quad (3)$$

To solve for R, we first notice that since $h(X_n|X_{n+1}) = 0$ (assumption 2) we can show by induction that

$$h(X_{n+1}|X_n, X_{n-1}, \dots, X_1, X_0) = h(X_{n+1}|X_n) \quad (4)$$

Due to the symmetric nature of mutual information, we can prove the equation below [2].

$$h(X_{n+1}|X_n) = h(X_{n+1}) - h(X_n) + h(X_n|X_{n+1}) \quad (5)$$

Bringing us to the equation for R below

$$R = \lim_{n \to \infty} \frac{h(X_{n+1}) - h(X_n)}{dt} \quad (6)$$

Next, we use assumption 3 to re-write the difference in entropy at time step n and n-1 as equal to the difference in entropy at time step 1 and the initial state.

$$R = \lim_{n \to \infty} \frac{h(X_1) - h(X_0)}{dt} \quad (7)$$

Since the $X_n$'s are Gaussian, we can easily calculate the differential entropy of each step using equation (2) and the differential entropy of the Gaussian distribution [2].

$$R = \lim_{n \to \infty} \left[\frac{\frac{1}{2}\log(2\pi e \cdot (\Delta x_0 + \Delta p \cdot dt/m)^2) - \frac{1}{2}\log(2\pi e \cdot (\Delta x_0)^2)}{dt}\right] (8)$$

$$R = \lim_{n \to \infty} \left[\log\left(1 + \frac{\Delta p \cdot dt}{\Delta x_0 m}\right)\right] \quad (9)$$

Using equations (31) and (32) this is re-written

$$R = \lim_{n \to \infty} \left[\log\left(1 + \frac{2k_B T \cdot dt}{\hbar}\right)\right] \quad (10)$$

We are assured by assumption 1 that $dt \ll \hbar/(2k_B T)$. Thus, we can Taylor expand the logarithm giving the first term plus the terms that are O(dt) or smaller.

$$R = \lim_{n \to \infty} \left[\frac{2k_B T}{\hbar} + O(dt) + \cdots\right] (11)$$

Ignoring the terms of O(dt) or smaller, we get our primary result

$$R = 2k_B T/\hbar \quad (12)$$

The other method to calculate the entropy rate is $R'$, which equals the limit as n goes to infinity of the entropy of all the $X_n$'s divided by n times dt [3]. Since we are looking at the rate of generation of the entropy (not the initial conditions), we subtract the entropy of the initial state $h(X_0)$. This also assures that R is in the correct units.

$$R' \equiv \lim_{n \to \infty} \left[\frac{h(X_n, X_{n-1}, \dots, X_1, X_0) - h(X_0)}{n \cdot dt}\right] \quad (13)$$

Since $h(X_n|X_{n+1}) = 0$ (assumption 2), we know that $h(X_n, X_{n-1}, \dots, X_1, X_0) = h(X_n)$, thus

$$R' = \lim_{n \to \infty} \left[\frac{h(X_n) - h(X_0)}{n \cdot dt}\right] \quad (14)$$

We can insert zero into the limit

$$0 = \sum_{k=1}^{n-1} -h(X_k) + h(X_k) \quad (15)$$

$R'$ becomes

$$R' = \lim_{n \to \infty} \left[\frac{\sum_{k=1}^{n} h(X_k) - h(X_{k-1})}{n \cdot dt}\right] (16)$$

Assumption 3 now lets us rewrite this as

$$R' = \lim_{n \to \infty} \left[\frac{n \cdot (h(X_1) - h(X_0))}{n \cdot dt}\right] (17)$$

We see that

$$R = R' \quad (18)$$

We can safely conclude that

$$R = R' = 2k_B T/\hbar \quad (19)$$

In this view, the temperature acts as an average energy and generates information (or entropy) at a rate equal to twice the average energy divided by ℏ.

## 5. The Variance of $X_n$:

Given the wave particle duality, which states that a free particle is both a wave and a particle, we see that our free particle undergoes both quantum mechanical diffusion of the wave and classical diffusion of the particle.

Introducing $X_p$, $X_f$, $p(x,t)$ and $f(x,t)$ makes this more clear. $X_p$ is a random variable drawn from $p(x,t) = \psi^*(x,t)\psi(x,t)$, the probability distribution associated with the quantum mechanical wave-function, which is the solution to the quantum diffusion equation, equation (33). $X_f$ is a random variable drawn from $f(x,t)$ and is the solution to real diffusion equation, equation (42).

If $X_g$ were an observation of where the particle is located, it would be the sum of a sample $X_p$ drawn from $p(x,t)$ and the uncorrelated sample $X_f$, drawn from $f(x,t)$.

$$X_g = X_p + X_f \quad (20)$$

Thus the action of $f(x,t)$ is to translate the center of the wave function, $\psi(x,t)$, by a sample of $X_f$.

As we know from probability theory, the resulting distribution, $g(x,t)$ is equal to the convolution of $p(x,t)$ and $f(x,t)$ over the x variable (30) [4].

$$g(x,t) = p(x,t) *_x f(x,t) \quad (21)$$

Since both $p(x,t)$ and $f(x,t)$ are Gaussian distributions, it is easy to show that the convolution of the two is again a Gaussian distribution with an expected value being equal to the sum of the two expected values (which in this case is zero) and a variance that is equal to the sum of the variances of the individual distributions.

$$\overline{x_g} = \overline{x_p} + \overline{x_f} = 0 \quad (22)$$
$$(\Delta x_g)^2 = (\Delta x_p)^2 + (\Delta x_f)^2 \quad (23)$$

Shown in equation (40) the variance of $p(x,t)$ is $(\Delta x_p)^2$.

$$(\Delta x_p)^2 = (\Delta x_0)^2 + \left(\Delta p \cdot \frac{t}{m}\right)^2 \quad (24)$$

In this equation t is the amount of time that has passed since the particle was initialized in the minimum uncertainty state, $\Delta x_0$ is the standard deviation of the minimum uncertainty state, $\Delta p$ is the standard deviation of the minimum uncertainty state in the momentum domain and m is the mass of the particle.

Shown in equation (48), the variance of $f(x,t)$ is $(\Delta x_f)^2$.

$$(\Delta x_f)^2 = \frac{\hbar \cdot t}{m} \quad (25)$$

Thus we get $(\Delta x_g)^2$.

$$(\Delta x_g)^2 = (\Delta x_0)^2 + \left(\Delta p \cdot \frac{t}{m}\right)^2 + \frac{\hbar \cdot t}{m} \quad (26)$$

Inserting into the last term the Heisenberg Uncertainty principle (32), $2\Delta x_0 \Delta p = \hbar$, we can group.

$$(\Delta x_g)^2 = \left(\Delta x_0 + \Delta p \cdot \frac{t}{m}\right)^2 \quad (27)$$

To understand the model, it is helpful to look at equation (26). $(\Delta x_g)^2$ is the sum of three variances. The first is from the Heisenberg Uncertainty Principle of the initialized state, the second is from the thermal drift of the center of the minimum uncertainty wave packet moving with a group momentum taken as a sample of the momentum domain, and the third is from the classical diffusion of the center of the wave-function on top of the other two.

It is also possible to derive equation (27) by assuming no force on the particle, which lets you deduce $x = x_0 + p \cdot t/m$. Squaring and taking the ensemble average is all you need [5].

## 6. The Imaginary Diffusion Equation:

The Kinetic Energy Hamiltonian characterizes the wave packet of a free particle in one dimension, where H is the Hamiltonian, p is the momentum along the x direction, and m is the mass of the particle [6].

$$H = \frac{p^2}{2m} \quad (28)$$

Given that the momentum commutes with the Hamiltonian, $[p,H] = [p,p^2/2m] = 0$, each eigenvalue of the momentum is a constant of motion and thus the variance in momentum space does not grow with time. It is possible to learn the width of the variance of the momentum by looking at the equipartition of energy [7]. Using the equipartition of energy we know to equate the degree of freedom associated with the average Kinetic Energy to one half the temperature times Boltzmann's constant.

$$\overline{\frac{p^2}{2m}} = \tfrac{1}{2} k_B T \quad (29)$$

Since we will assume that the average momentum is zero, we can solve for the variance of the momentum.

$$\bar{p} = 0 \quad (30)$$
$$\Delta p^2 = \overline{p^2} - \bar{p}^2 = k_B T m \quad (31)$$

Also from the Heisenberg Uncertainty Principal, we can solve for the standard deviation of the wave-function in the spatial domain in terms of its width in the momentum space.

$$\Delta x_0 = \frac{\hbar}{2\Delta p} \quad (32)$$

With these dependencies stated, we can move onto the imaginary diffusion equation, which takes the original Hamiltonian and rewrites it in terms of operators. Interpreting the Hamiltonian as the imaginary time derivative operator and the momentum as the negative imaginary spatial derivative operator we can take equation (28) and arrive at the imaginary diffusion equation

$$i\hbar \frac{d}{dt}\psi(x,t) = \frac{-\hbar^2}{2m}\frac{d^2}{dx^2}\psi(x,t) \quad (33)$$

Don't forget that we still have the eigenvalue equations (34,35) where H and p are the operators and ω and k are the eigenvalues.

$$H\psi(x,t) = \hbar\omega\psi(x,t) \quad (34)$$
$$\frac{p^2}{2m}\psi(x,t) = \frac{(\hbar k)^2}{2m}\psi(x,t) \quad (35)$$

We can calculate the different eigenvalues, $\omega$ and $k$, through equations (28,34,35) and as we should expect arrive at the equation for kinetic energy.

$$\hbar\omega = \frac{(\hbar k)^2}{2m} \quad (36)$$

To solve equation (33), we will begin in the momentum domain $\Psi(k/2\pi)$ and take the inverse Fourier Transform to observe how $\psi(x,t)$ evolves over time [8]. We use $k/2\pi$ (the wavenumber divided by $2\pi$) as the independent variable because we want both $\Psi(k/2\pi)$ and $\psi(x,t)$ to be normalizable to one.

$$\Psi(k/2\pi) = \left(\frac{2\pi}{(\Delta k)^2}\right)^{\frac{1}{4}} \cdot \exp\left(\frac{-k^2}{4(\Delta k)^2}\right) \quad (37)$$

Our assumption that the wave-function of the free particle in the momentum space is a Gaussian wave-packet is quite reasonable given the nice properties of the Gaussian. Similarly, this assumption is already implicit in the equipartition of energy which was used to find the width of the initial wave-packet. Because the equipartition theorem is derived from the perfect gas law (where particles are modeled using the binomial distribution, of which the Gaussian is the limit), the Gaussian is the right distribution to start with.

To properly account for the evolution of $\psi(x,t)$ governed by equation (33), $\exp(i(kx - \omega t))$ is used as the kernel for the inverse Fourier Transform.

$$\psi(x,t) = \int_{-\infty}^{\infty}\left(\frac{2\pi}{(\Delta k)^2}\right)^{\frac{1}{4}}\exp\left(\frac{-k^2}{4(\Delta k)^2}\right)\exp(i(kx - \omega t))\frac{dk}{2\pi} \quad (38)$$

Using equation (36) to substitute in for ω you can solve for equation (38) by completing the squares to get $\psi(x,t)$ [8]. $\psi(x,t)$ is in Gaussian form; to calculate the variance, we need to take the magnitude squared of the wave-function and get the distribution of the particle.

$$p(x,t) = \psi^*(x,t)\psi(x,t) = \frac{1}{\sqrt{2\pi(\Delta x_p)^2}}\exp\left(\frac{-x^2}{2(\Delta x_p)^2}\right) \quad (39)$$

Where

$$(\Delta x_p)^2 = (\Delta x_0)^2 + \left(\Delta p \frac{t}{m}\right)^2 \quad (40)$$

And

$$\Delta p = \hbar \Delta k \quad (41)$$

This is, of course, the well know result from quantum mechanics where the variance of the particle is the sum of the initial variance from the Heisenberg Uncertainty Principal and the associated variance of the momentum domain imparting a thermal group velocity $\Delta p/m$ [9].

## 7. The Real Diffusion Equation:

When the diffusion constant of a diffusion process is real and does not vary with position, the resulting diffusion equation is as below [10].

$$\frac{d}{dt}f(x,t) = D\frac{d^2}{dx^2}f(x,t) \quad (42)$$

Of course the solution to this real diffusion equation is the Gaussian with variance equal to 2Dt [11].

$$f(x,t) = \frac{1}{\sqrt{4\pi Dt}}exp\left(\frac{-x^2}{4Dt}\right) \quad (43)$$
$$(\Delta x_f)^2 = 2Dt \quad (44)$$

To find D, we will start with the imaginary diffusion operator and using analytical continuation, perform a Minkowski transformation [6]. The imaginary diffusion operator (33) is

$$\frac{d}{dt} = \frac{i\hbar}{2m}\frac{d^2}{dx^2} \quad (45)$$

Upon applying the Minkowski transformation, imaginary time is replaced with real time, $i \cdot t \to t'$ [12]. Applied on the imaginary diffusion operator, the Minkowski transformation brings out the real diffusion constant we are looking for.

$$\frac{d}{dt} = \frac{i\hbar}{2m}\frac{d^2}{dx^2} \to \frac{d}{dt'} = \frac{\hbar}{2m}\frac{d^2}{dx^2} \quad (46)$$

By observation we see that

$$D = \frac{\hbar}{2m} \quad (47)$$

We can also derive $D$ from kinematic arguments as was shown in [1]. We can calculate the variance of f(x,t).

$$(\Delta x_f)^2 = \frac{\hbar t}{2m} \quad (48)$$

## 8. Entropy at $t = 0$

It is important to ask the entropy of the initial state. We find that at $t = 0$ the entropy is 1 natural unit. Since the wavefunction and associated probability distribution are continuous, we calculate the entropy using the equation for differential entropy. One might object that the differential entropy is only accurate up to a scale factor. However I argue (and so did Hirshman [13]) that if you add the differential entropy in the dual domain, the scale factor cancels out because of the scale property of the Fourier Transform and the result is an absolute measure. As before we will use the position and the wavenumber divided by $2\pi$ as the dual domains.

$$h(\psi^*(x,t)\psi(x,t)) + h(\Psi^*(k/2\pi)\Psi(k/2\pi)) \quad (49)$$

Hirschman [12] showed that this entropy for any wave-function is $\geq \ln(e/2)$ and $= \ln(e/2)$ when the wave-function is Gaussian.

While we are working with a Gaussian initial state, the answer appears to be a little more complex than just $= \ln(e/2)$. We learn from [1] that when solving for the quantum and relativistic length scales of dark particles, particles come in pairs. With only one particle and no reference frame there is no way of knowing the position or the momentum, even if there were universal measuring sticks.

We get around this with two particles and a measuring stick/clock by determining the relative displacement and speed. Thus we need to look at the entropy of the relative difference of position and momentum of the two particles.

Define $p_1 = p(x_1, t=0) = \psi^*(x_1, 0)\psi(x_1, 0)$ as the probability distribution on the location of particle 1 at $t = 0$ and similarly for $p_2$ for particle 2. For the momentum space define $P_1 = \Psi^*(k_1/2\pi)\Psi(k_1/2\pi)$ as the probability distribution on the wavenumber divided by $2\pi$ for the first particle and $P_2$ for the second particle.

The probability distribution on the relative displacement and wavenumber, $\delta_x = x_1 - x_2$ and $\delta_k = k_1 - k_2$ are $p_{\delta_x}$ and $P_{\delta_k}$, respectively, and will be Gaussian assuming both the reference particle and the initial particle have Gaussian wave-functions. Since differential entropy is invariant to the first moment we can assume without loss of generality, the first moment of the reference wave-function is zero. The second moment of $\delta_x$ and $\delta_p$ will be the sum of the respective second moments of the particle and reference particle if the two are not correlated.

We can go even further and show that the reference particle should have the same second moments as the particle we are measuring if we minimize the entropy. Thus, we arrive at the distributions for both domains for $\delta_x$ and $\delta_p$

$$p_{\delta_x}(\delta_x)d\delta_x = \frac{1}{\sqrt{4\pi(\Delta x_0)^2}} exp\left(\frac{-\delta_x^2}{4(\Delta x_0)^2}\right)d\delta_x \quad (49)$$

$$P_{\delta_k}(\delta_p)d\delta_p = \frac{1}{\sqrt{4\pi(\Delta k/2\pi)^2}} exp\left(\frac{-(\delta_p/2\pi)^2}{4(\Delta k/2\pi)^2}\right)\frac{d\delta_p}{2\pi} \quad (50)$$

Thus, the total absolute entropy of the initial state, $I_0$, is

$$I_0 = h\left(p_{\delta_x}(\delta_x)\right) + h\left(P_{\delta_k}(\delta_p/2\pi)\right) \quad (51)$$

$$I_0 = \frac{1}{2}\log(2\pi e \cdot 2(\Delta x_0)^2) + \frac{1}{2}\log(2\pi e \cdot 2(\Delta k/2\pi)^2) \quad (52)$$

$$I_0 = 1 \; natural \; unit \quad (53)$$

There are 2 things of note relative to the rate, $R$, calculated above. First we see that since the rate, $R$, from above, is the difference between the entropy at two times, the impact of the wider distribution of $\delta_x$ vs. $x$ is negated. Thus we could have done the analysis above using $\delta_x$ and $\delta_p$ instead $x$ and $p$ and the result would be the same. Second, we see that the entropy of the initial state is equal to the additional entropy generated by the diffusion process during the de-coherence time, $\tau = \hbar/2k_B T$

### 9. Conclusion:

We have seen that by making three assumptions about the thermal diffusion of a free particle, we are able to show that entropy is generated at a rate equal to twice the particle's temperature (when expressed in the correct units).

This result will be applicable to all studies on free particles and other environments that are governed by similar equations. Also a myriad of applications exist in computer modeling, including but not limited to the following: finite difference time domain methods, Block's equations for nuclear magnetic resonance imaging, and plasma and semiconductor physics.

To check the primary result, one would perform a quantum non-demolition measurement on the quantum state of an ensemble of free particles. The minimum bit rate needed to describe the resulting string of numbers that describe the trajectory would be the entropy rate and should be equal to twice the temperature.

However even before an experiment can be conducted, this result is useful by suggesting the use of different information theoretical techniques to examine problems with de-coherence and might give a different perspective on the meaning of temperature.

This result is interesting as a stand-alone data point, that the entropy rate is equal to twice the temperature. However if we could to go further and more generally say that temperature is the same as entropy rate, it would change the way we view temperature and entropy.

Kitzbühel, Austria 2008

### 10. Acknowledgements

JLH thanks Thomas Cover for sharing his passion for the Elements of Information Theory, and McKinsey & Co. for the amazing environment where this article was written.

### REFERENCES

1. J. Haller Jr., "Dark Particles Answer Dark Energy," *Journal of Modern Physics*, Vol. 4 No. 7A1, 2013, pp. 85-95. doi: 10.4236/jmp.2013.47A1010
2. Yamamoto & Imamoglu, Mesoscopic Quantum Optics, John Wiley & Sons, New York, 1999
3. Cover& Thomas, Elements of Information Theory, John Wiley & Sons, New York 1991
4. Bracewell, The Fourier Transform and Its Applications, 2nd ed.,McGraw Hill, New York 1986
5. Gardiner, C. W. & Zoller, P., Quantum Noise, Springer, Berlin, 2004
6. Shankar, Principles of Quantum Mechanics, Plenum Press, New York 1994
7. Feynman, Lectures on Physics, Addison-Wesley Publishing, Reading Massachusetts, 1965

8. Bohm, <u>Quantum Theory</u>, Dover Publications, Mineola, N.Y. 1989
9. Shankar, <u>Principles of Quantum Mechanics</u>, Plenum Press, New York 1994
10. Bittencourt, <u>Fundamentals of Plasma Physics, 2nd ed.</u>, Sao Jose dos Campos, SP 1995
11. Einstein, "Investigation on the Theory of, The Brownian Movement", Translated by Cowper, Dover 1956
12. Einstein, <u>The Meaning of Relativity, 5th Edition</u>, Princeton University Press, Princeton NJ 1956
13. I. I. Hirshman, "A Note on Entropy," *Amer. J. of Math*. January 1957, 79, No 1, p. 152